\documentclass[12pt,titlepage]{article}
\usepackage{amssymb}
\input epsf
\newcommand{\proof}{{\noindent \bf Proof. }}
\newtheorem{thm}{Theorem}
\newtheorem{defi}{Definition}

\newtheorem{cor}[thm]{Corollary}

\newcommand{\F}{{\cal F}}

\def\2{\mathbb Z_2}
\def\ind{{\rm ind}}
\def\coind{{\rm coind}}

\def\qed{$\Box$}

\def\cd{{\rm cd}}
\def\KG{{\rm KG}}
\def\SG{{\rm SG}}
\newcommand\mbf[1]{\mbox{\boldmath$#1$}}
\newcommand\msbf[1]{\mbox{\boldmath\scriptsize$#1$}}

\title{Colorful subgraphs in Kneser-like graphs}
\author{
      {\bf G\'abor Simonyi}\thanks{Research partially supported by the
Hungarian Foundation for Scientific Research Grant (OTKA) Nos.\ T037846,
T046376, and AT048826.}\\
 Alfr\'ed R\'enyi Institute of Mathematics\\
 Hungarian Academy of Sciences\\
 1364 Budapest, POB 127, Hungary\\
 {\tt simonyi@renyi.hu}
\and
   {\bf G\'abor Tardos}\thanks{Research partially supported by the
Hungarian Foundation for Scientific Research Grant (OTKA) Nos.\ T037846,
T046234, and AT048826.}\\
 School of Computing Science \\
 Simon Fraser University \\
 Burnaby BC, Canada V5A 1S6\\
 and\\
 Alfr\'ed R\'enyi Institute of Mathematics\\
 Hungarian Academy of Sciences\\
 1364 Budapest, POB 127, Hungary\\
 {\tt tardos@cs.sfu.ca}}

\begin{document}
\maketitle

\begin{abstract}
Combining Ky Fan's theorem with ideas of Greene and Matou\v{s}ek we prove a
generalization of Dol'nikov's theorem. Using another variant of the
Borsuk-Ulam theorem due to Bacon and Tucker, we also prove the presence of all
possible completely multicolored $t$-vertex complete bipartite graphs in
$t$-colored $t$-chromatic Kneser graphs and in several of their relatives.
In particular, this implies a generalization of a recent result of G. Spencer
and F.\ E. Su.

\end{abstract}

\section{Introduction}

The solution of Kneser's conjecture in 1978 by L\'aszl\'o Lov\'asz \cite{LLKn}
opened up a new area of combinatorics that is usually referred to as the
topological method \cite{Bjhand}. Many results of this area, including the
first one by Lov\'asz, belong to one of its by now most developed branches
that applies the celebrated Borsuk-Ulam theorem to graph coloring problems. An
account of such results and other applications of the Borsuk-Ulam theorem in
combinatorics is given in the excellent book of Matou\v{s}ek \cite{Mat}.

Recently it turned out that a generalization of the Borsuk-Ulam theorem found
by Ky Fan \cite{kyfan} in 1952 can give useful generalizations and variants of
the Lov\'asz-Kneser theorem. Examples of such results can be found in
\cite{Meunier, ST, Spthes, SpSu}.

In this note our aim is twofold. In Section~\ref{gdol} we show a further
application of Ky Fan's theorem. More precisely, we give a generalization
of Dol'nikov's theorem, which is itself a generalization of the Lov\'asz-Kneser
theorem. The proof will be a simple combination of Ky Fan's result with
the simple proof of Dol'nikov's theorem given by Matou\v{s}ek in \cite{Mat}
that was inspired by Greene's recent proof \cite{Gre} of the Lov\'asz-Kneser
theorem.

In Section~\ref{BT} we use another variant of the Borsuk-Ulam theorem due to
Philip Bacon \cite{Bacon} and A. W. Tucker \cite{Tuc} to show a property of
optimal
colorings of certain $t$-chromatic graphs, including Kneser graphs, Schrijver
graphs , Mycielski graphs, Borsuk graphs, odd chromatic rational complete
graphs, and two other types of graphs appearing in \cite{ST} that will also be
defined in Subsection~\ref{gr}.
The claimed property (in a somewhat weakened form) is that all the complete
bipartite graphs $K_{l,m}$ with $l+m=t$ will have totally multicolored copies
in proper $t$-colorings of the above graphs. 

When applied to rational complete graphs this implies a new proof of an
earlier result about the circular chromatic number (see
Subsection~\ref{rat}). It was originally obtained in \cite{ST} and
independently in \cite{Meunier} for some special cases.  

When applied to Kneser graphs the above property 
implies a generalization of a recent result due to Gwen Spencer and Francis
Edward Su \cite{Spthes, SpSu} which we will also present in the last
subsection.

\section{A generalization of Dol'nikov's theorem} \label{gdol}

We recall some concepts and notations from \cite{Mat}.
For any family $\F$ of subsets of a fixed finite set we define the general
Kneser graph $\KG(\F)$ by
\begin{eqnarray*}
V(\KG(\F))&=&\F,\\
E(\KG(\F))&=&\{\{F,F'\}: F,F'\in\F, F\cap F'=\emptyset\}.
\end{eqnarray*}
We note about our terminology that
when referring to a Kneser graph (without the adjective ``general'') we mean
the ``usual'' Kneser graph $\KG(n,k)$ that is identical to the general Kneser
graph of a set system consisting of all $k$-subsets of an $n$-set.

A hypergraph $H$ is $m$-colorable, if its vertices can be colored by (at most)
$m$ colors so that no hyperedge becomes monochromatic. The $m$-colorability
defect of a set system $\emptyset\notin\F$ (identified with a hypergraph in
the obvious way) is
$$\cd_m(\F):=\min\left\{|Y|: (X\setminus Y, \{F\in \F: F\cap Y=\emptyset\})\
{\rm is}\ m{\rm -colorable}\right\}.$$

\medskip
\par
\noindent
{\bf Dol'nikov's theorem} (\cite{Dol}) {\it For any finite set system
  $\emptyset\notin\F$, 
  the inequality $$\cd_2(\F)\leq \chi(\KG(\F))$$ holds.}

\bigskip

This theorem generalizes the Lov\'asz-Kneser theorem, as it is easy to check
that if $\F$ consists of all the $k$-subsets
of an $n$-set with $k<n/2$, then $\cd_2(\F)=n-2k+2$ which is the true value of
$\chi(\KG(\F))$ in this case.
On the other hand, as also noted in \cite{Mat}, equality
between $\cd_2(\F)$ and $\chi(\KG(\F))$ does not hold in general.

Recently Greene \cite{Gre} found a very simple new proof of the
Lov\'asz-Kneser theorem. (This proof
follows the footsteps of B\'ar\'any's \cite{Bar} earlier simple proof but with
a skillful trick it can avoid the use of Gale's theorem.) In \cite{Mat}
Matou\v{s}ek observed that one can generalize Greene's proof so that it also
gives Dol'nikov's theorem. Here we combine this proof with Ky Fan's theorem
(see below) to obtain the following generalization.

\begin{thm} \label{thm:gdol}
Let $\F$ be a finite family of sets, $\emptyset\notin
\F$ and $\KG(\F)$ its
general Kneser graph. Let $r=\cd_2(\F)$. Then any proper coloring of $\KG(\F)$
with colors
$1,\dots,m$ ($m$ arbitrary) must contain a completely multicolored complete bipartite
graph $K_{\lceil r/2\rceil,\lfloor r/2\rfloor}$ such that the $r$ different
colors occur alternating on the two sides of the bipartite graph with respect
to their natural order.
\end{thm}

This theorem generalizes Dol'nikov's theorem, because it implies that any
proper coloring must use at least $\cd_2(\F)$ different colors.

\medskip
\par
\noindent {\it Remark 1.}\ Theorem~\ref{thm:gdol} is clearly in
the spirit of the Zig-zag theorem of \cite{ST} the special case
of which for Kneser graphs was already established by Ky Fan in
\cite{kyfan2}. This theorem claims that if a specific topological
parameter of a graph $G$, the value of which is a lower bound on
its chromatic number, is at least $t$ then any proper coloring of
the graph must contain a completely multicolored $K_{\lceil
t/2\rceil,\lfloor t/2\rfloor}$ where the colors also alternate on
the two sides with respect to their natural order. As we will show
in Remark~2 the proof below can be modified to show that the
topological parameter mentioned is at least $\cd_2(\F)$ for any
$\KG(\F)$. We will say more about this topological parameter in
Section~\ref{BT}. $\hfill\Diamond$

\medskip

To prove the theorem we first have to state Ky Fan's theorem
\cite{kyfan}. Just like the Borsuk-Ulam theorem it
has several equivalent forms (see \cite{kyfan, Bacon}). The one
fitting best for our purposes is the following. (To find exactly this form,
see \cite{Bacon}.) For a set $A$ on the unit sphere $S^h$ we denote by $-A$ its
antipodal set, i.e., $-A=\{-\mbf{x}: \mbf{x}\in A\}$.

\medskip
\par
\noindent
{\bf Ky Fan's theorem} (\cite{kyfan}) {\it Let $A_1,\dots,A_m$ be open
  subsets of
  the $h$-dimensional sphere $S^h$ satisfying that none of them contains
  antipodal points (i.e., $\forall i$ and $\forall
  \mbf{x}\in S^h\ \mbf{x}\in A_i$ implies $-\mbf{x}\notin A_i$) and that at
  least one of $\mbf{x}$ and $-\mbf{x}$ is contained in $\cup_{i=1}^m A_i$ for
  all $\mbf{x}\in S^h$.

Then there exists an $\mbf{x}\in S^h$ and $h+1$ distinct indices $i_1<\dots
<i_{h+1}$ such that $\mbf{x}\in A_{i_1}\cap -A_{i_2}\cap\dots\cap(-1)^h
A_{i_{h+1}}$.}

\bigskip
\par
\noindent
{\bf Proof of Theorem~\ref{thm:gdol}.}
Let $h=\cd_2(\F)-1$ and consider the sphere $S^h$. We assume without loss
of generality that the base set $X:=\cup\F$ is finite and identify its
elements with points of $S^h$ in
general position, i.e., so that at most $h$ of them can be on a common
hyperplane through the origin. Consider an arbitrary fixed proper
coloring of $\KG(\F)$ with colors $1,\dots,m$.
For every $\mbf{x}\in S^h$ let $H(\mbf{x})$ denote the open hemisphere
centered at $\mbf{x}$.
Define the sets $A_1,\dots,A_m$ as follows. Set
$A_i$ will contain exactly those points $\mbf{x}\in S^h$ that have the
property that $H(\mbf{x})$
contains the points of some $F\in\F$ which is colored by color $i$ in the
coloring of $\KG(\F)$ considered. The sets $A_i$ are all open. None of them
contains an antipodal pair of points, otherwise there would be two disjoint
open hemispheres both of
which contain some element of $\F$ that is colored $i$. But these two elements
of $\F$ would be disjoint contradicting the assumption that the coloring was
proper. Now we show that there is no $\mbf{x}\in S^h$ for which neither
$\mbf{x}$ nor $-\mbf{x}$ is in $\cup_{i=1}^m A_i$. Color the points in $X\cap
H(\mbf{x})$ red, the points in $X\cap H(-\mbf{x})$ blue and delete the points
of $X$ not colored, i.e., those on the ``equator'' between $H(\mbf{x})$ and
$H(-\mbf{x})$. Since at most $h<\cd_2(\F)$ points are deleted there exists some
$F\in\F$ which became completely red or completely blue. All points of $F$
are either in $H(\mbf x)$ or in $H(-\mbf x)$. This implies that $\mbf x$ or
$-\mbf x$ should belong to $A_i$ where $i$ is the color of $F$ in our fixed
coloring of $\KG(\F)$.

Thus our sets $A_1,\dots,A_m$ satisfy the conditions and therefore also the
conclusion of Ky Fan's theorem. Let $F_{i_j}\in \F$ be the set responsible for
$\mbf{x}\in A_{i_j}$ for $j$ odd and for $-\mbf{x}\in A_{i_j}$ for
$j$ even in the conclusion of Ky Fan's theorem. Then all the $F_{i_j}$'s with
odd $j$ must be disjoint from all the $F_{i_j}$'s with even $j$. Thus they
form the complete bipartite graph claimed. \hfill\qed

\subsection{Consequences}

We would like to point out that Theorem~\ref{thm:gdol} can really be
stronger than Dol'nikov's theorem, especially if looked at as an upper
bound for $\cd_2(\F)$. It is not difficult to see (and is included as an
exercise in Matou\v{s}ek's book \cite{Mat}, cf.\ page 64, Exercise 3) that
every graph $G$ is isomorphic to the general Kneser graph $\KG(\F)$ of some set
system
$\F$.

Consider a graph $G$ with girth at least $5$ and high chromatic
number. (According to a famous result of Erd\H{o}s \cite{Erd} such graphs
exist.) Let $\F_G$ denote a set system with the property that $\KG(\F_G)=G$.
Theorem~\ref{thm:gdol} implies that $\cd_2(\F_G)\leq 3$ since $\KG(\F_G)$ does
not contain any $K_{2,2}$ subgraph. At the same time, the upper bound
Dol'nikov's
theorem implies on $\cd_2(\F_G)$ is the chromatic number of $G$ which may be
one million.

Formulating the above argument about $K_{2,2}$-free graphs in terms of the set
system, we have the following.

\begin{cor} \label{cd3}
Let $\F$ be a set system not containing the empty set and satisfying that for
any two distinct sets $A,B\in\F$ there is at most one set in $\F$ that is
disjoint from both $A$ and $B$. Then $\cd_2(\F)\leq 3.$
\end{cor}

\proof The condition on $\F$ implies that $K_{2,2}\nsubseteq \KG(\F)$, thus
Theorem~\ref{thm:gdol} implies the statement. \hfill\qed

\medskip

It would be unfair to deny that Corollary~\ref{cd3} can be proven easily in a
rather elementary way, too. Just consider sets $A,B\in \F$ for which $A\cup B$
is minimal. By this choice any $C\subseteq A\cup B$, $C\in \F$, should satisfy
$(A\setminus B)\cup (B\setminus A)\subseteq C$. Therefore deleting a point
from $A\setminus B$ and
one from $B\setminus A$ if neither of these is empty and coloring the rest of
$A\cup B$ red, no completely red element of $\F$ may occur. There is at most
one $D\in \F$ with $D\cap(A\cup B)=\emptyset$. Delete a point of $D$ (if it
exists) and color all the remaining points blue. This gives a red-blue
coloring that proves the statement in case $A\setminus B\neq\emptyset$ and
$B\setminus A\neq\emptyset$ both hold. If, say, $A\setminus B=\emptyset$, then
we have to delete a point of $A$ instead of one of $A\setminus B$ that is
impossible now. Otherwise we can do everything similarly thereby proving the
statement also in this other case.

Even though such an elementary proof exists we believe that the condition of
Corollary~\ref{cd3} comes to mind much more naturally because of the
topological background. Also, Theorem~\ref{thm:gdol} implies other statements
of the like, that are perhaps more complicated to prove in an elementary way.
For example, even if $\KG(\F)$ contains one $K_{q,q}$ subgraph for some huge
$q$, but no other $K_{2,2}$ subgraph (i.e., no $K_{2,2}$ subgraph apart from
those contained in the large $K_{q,q}$) then $\cd_2(\F)\leq 3$ still
holds. This is
because the colors appearing on the two sides of the $K_{q,q}$ can be labeled
so that the colors on one side all precede the colors of the other side in
their ordering and therefore a $K_{2,2}$ required by Theorem~\ref{thm:gdol}
for $\cd_2(\F)\ge 4$ still cannot occur.

\section{Applying a theorem of Bacon and Tucker} \label{BT}

\subsection{Preliminaries} \label{prel}

First we give a very brief introduction of some topological concepts we
need. All this can be found in detail, e.g., in \cite{Mat}.
A $\2$-space is a pair $(T,\nu)$, where $T$ is a topological space and
$\nu:T\to T$ is an involution, that is, a continuous map satisfying
$\nu(\nu(x))=x$ for all $x\in T$. A $\2$-space $(T,\nu)$ is free if
$\nu(x)\neq x$ for every $x\in T$. If $\nu$ is clear from the context we
write $T$ in place of $(T,\nu)$. Accordingly, we write $S^h$ for the most
important $\2$-space we deal with, the $h$-dimensional sphere with the usual
antipodal map as involution.

A continuous map $f: (T,\nu)\to (W,\mu)$ is a $\2$-map if it respects
the respective involutions, that is, $f(\nu(x))=\mu(f(x))$ for every $x\in
T$. We write $(T,\nu)\to (W,\mu)$ if there exists a $\2$-map from $(T,\nu)$ to
$(W,\mu)$.
Two important parameters of a $\2$-space are its $\2$-index and
$\2$-coindex that are defined as
$$\ind(T,\nu):=\min\{h\ge 0:(T,\nu)\to S^h\},$$
and
$$\coind(T,\nu):=\max\{h\ge 0:S^h\to(T,\nu)\},$$
respectively.
The inequality $$\coind(T,\nu)\leq \ind(T,\nu)$$ always holds and is one of the
standard forms of the celebrated Borsuk-Ulam theorem.

In applications of the topological method one often associates so-called box
complexes to graphs. These give rise to topological spaces the index and
coindex of which can serve to obtain lower bounds for the chromatic number of
the graph. Following ideas in earlier works by Alon, Frankl, Lov\'asz
\cite{AFL} and others, the paper \cite{MZ} introduces several box complexes
two of which we also define below.

\begin{defi}
The box complex $B(G)$ is a simplicial complex on the vertices
$V(G)\times\{1,2\}$. For subsets $S,T\subseteq V(G)$ the set $S\uplus
T:=S\times\{1\}\cup T\times\{2\}$ forms a simplex if and only if
$S\cap T=\emptyset$, the vertices in $S$ have at least one common
neighbor and the same is true for $T$, and the complete bipartite graph with
sides $S$ and $T$ is a subgraph of $G$. The $\mathbb{Z}_2$-map $S\uplus
T\mapsto T\uplus S$ acts
simplicially on $B(G)$ making the body $||B(G)||$ of the complex a free
$\mathbb{Z}_2$-space.
\end{defi}

It is explained in \cite{MZ} and \cite{Mat} that $B(G)$ is a functor, meaning
for example,
that whenever there exists a homomorphism from a graph $F$ to another graph
$G$ then $B(F)\to B(G)$ is also true. It is not hard to see that
$||B(K_n)||\cong S^{n-2}$ with a $\2$-homeomorphism (i.e., a homeomorphism
which is a $\2$-map). For the $\2$-index and $\2$-coindex of $||B(G)||$ we
simply write $\ind(B(G))$ and $\coind(B(G))$, respectively, and we will do
similarly for the other box complex $B_0(G)$ defined below.
Since a graph is $t$-colorable if and only if it admits
a homomorphism to $K_t$, the foregoing implies
\begin{equation}\label{eq:chib}
\chi(G)\ge \ind(B(G))+2\ge\coind(B(G))+2.
\end{equation}

Another box complex $B_0(G)$ defined in \cite{MZ} differs from
$B(G)$ only by containing all those simplices $S\uplus T$, too,
where one of $S$ or $T$ is empty independently of the existence
of common neighbors required in the definition of $B(G)$.

\begin{defi}
The box complex $B_0(G)$ is a simplicial complex on the vertices
$V(G)\times\{1,2\}$. For subsets $S,T\subseteq V(G)$ the set $S\uplus
T:=S\times\{1\}\cup T\times\{2\}$ forms a simplex if and only if
$S\cap T=\emptyset$, and the complete bipartite graph with sides $S$ and $T$ is
a subgraph of $G$. The $\mathbb{Z}_2$-map $S\uplus T\mapsto T\uplus S$ acts
simplicially on $B_0(G)$ making the body $||B_0(G)||$ of the complex a free
$\mathbb{Z}_2$-space.
\end{defi}

\medskip

Csorba \cite{Cs} proved a strong topological relationship between $B(G)$
and $B_0(G)$, namely, that the body of $B_0(G)$ is $\2$-homotopy
equivalent to the suspension of the body of $B(G)$. This extends
(\ref{eq:chib}) to the following longer chain of inequalities, cf.\ \cite{MZ}
and also \cite{ST}.
\begin{equation}\label{eq:chibb0}
\chi(G)\ge\ind(B(G))+2\ge\ind(B_0(G))+1\ge\coind(B_0(G))+1\ge\coind(B(G))+2
\end{equation}

Note that $B_0(G)$ is also a functor, and it is easy to see even without
Csorba's result that $||B_0(K_n)||\cong S^{n-1}$ with a $\2$-homeomorphism.

There are several interesting graph families the members of which satisfy the
inequalities in (\ref{eq:chibb0}) with equality. These include, for example,
Kneser graphs, and a longer list is given in Corollary~\ref{list} below.
(We note that some of the graphs in Corollary~\ref{list} give equality only in
the first three of the above inequalities, cf.\ Subsection~\ref{Umr}.)
\medskip

\medskip
\par
\noindent
{\it Remark 2.} The topological parameter mentioned in Remark 1
in Section~\ref{gdol} is $\coind(B_0(G))+1$. Thus our claim in Remark 1 is
that the proof of Theorem~\ref{thm:gdol} implies $\coind(B_0(\KG(\F)))\ge
\cd_2(\F)-1$ for any $\F$ not containing the empty set. Here we sketch the
proof of
this claim which is similar to the proof of Proposition 8 in \cite{ST}. Assume
again without loss of generality that $X=\cup\F$ is finite and identify its
elements with points of $S^h$ in general position as in the proof of
Theorem~\ref{thm:gdol} with $h=\cd_2(\F)-1$. For each vertex $v$ of $\KG(\F)$
and
$\mbf x\in S^h$ let $D_v(\mbf x)$ be the smallest distance of a point in $v$
(this point is an element of $X$) from the set $S^h\setminus H(\mbf
x)$. Notice that $D_v(\mbf x)>0$ iff $H(\mbf x)$ contains all points of
$v$. Set $D(\mbf x):=\sum_{v\in \F} (D_v(\mbf x)+D_v(-\mbf x))$. The argument
in the proof of
Theorem~\ref{thm:gdol} implies $D(\mbf x)>0$. Therefore the map $f(\mbf
x)=(1/D(\mbf x))(\sum_v D_v(\mbf x)||(v,1)||+\sum_v D_v(-\mbf x)||(v,2)||)$ is
a $\2$-map from $S^h$ to $||B_0(\KG(\F))||$, thus $\coind(B_0(\KG(\F)))\ge h$
as claimed.
\hfill$\Diamond$

\subsection{A colorful $K_{l.m}$-theorem} \label{klm}

In their recent paper \cite{CsLSW} Csorba, Lange, Schurr, and Wassmer proved
that if $\ind(B(G))\ge l+m-2$ then $G$ must contain the complete bipartite
graph $K_{l,m}$ as a subgraph. They called this ``the $K_{l,m}$-theorem''. In
case of those graphs that satisfy $\coind(B_0(G))+1=\chi(G)$ (for a list of
such graphs see Corollary~\ref{list}), the following statement generalizes
their result. We use the notation $[t]:=\{1,\dots,t\}$.

\begin{thm} \label{cklm}
Let $G$ be a graph for which $\chi(G)=\coind(B_0(G))+1=t$. Let $c: V(G)\to [t]$
  be a proper coloring of $G$ and let
  $A,B\subseteq [t]$ form a bipartition of the color set, i.e., $A\cup B=[t]$
  and $A\cap B=\emptyset$.

Then there exists a complete bipartite subgraph
$K_{l,m}$ of $G$ with sides $L,M$ such that $|L|=l=|A|,\ |M|=m=|B|$, and
$\{c(v):v\in L\}=A$, and $\{c(v):v\in M\}=B$. In particular, all vertices
of this $K_{l,m}$ receive different colors at $c$.
\end{thm}

For the proof we will use a modified version of the following theorem.
It is given in \cite{Bacon} for more general $\2$-spaces in place of $S^h$ but
for our purposes this restricted version will be sufficient.

\medskip
\par
\noindent
{\bf Bacon-Tucker theorem} (\cite{Bacon,Tuc})
{\it If $C_1,\dots,C_{h+2}$ are closed subsets of $S^h$, $$\cup_{i=1}^{h+2}
C_i=S^h, \ \ \forall i:\ C_i\cap(-C_i)=\emptyset,$$
and $j\in\{1,\dots,h+1\}$, then there is an $\mbf{x}\in S^h$ such that
$$\mbf{x}\in \cap_{i=1}^j C_i,\ {\rm and}\ -\mbf{x}\in\cap_{i=j+1}^{h+2}
C_i.$$}
\bigskip
\par

Bacon \cite{Bacon} shows the above theorem to be equivalent to $14$ other
statements that include standard forms of the Borsuk-Ulam theorem, and also
Ky Fan's theorem. Denoting its statement by $A_h(T)$ for a more general
$\2$-space $T$, he writes the following about the origins of this result: A
weakened form of $A_h(S^h)$ was stated in 1935 (\cite{AH}, Satz X, p. 487);
$A_2(S^2)$ is stated and the higher dimensional cases hinted at by Tucker in
1945 \cite{Tuc}. (The weakened form in \cite{AH} states only that any $h+1$ of
the $C_i$'s intersect which would be too weak for our purposes.)

It is a routine matter to see that the Bacon-Tucker theorem also holds for
open sets $C_i$. (For a collection of open sets $C_i$ one can define closed
sets $C_i'$ so that $C_i'\subseteq C_i$ for all $i$ and $\cup_{i=1}^{h+2} C_i'=
\cup_{i=1}^{h+2} C_i=S^h$. Cf.\ \cite{AH}, Satz VII, p. 73, quoted also in
\cite{Bacon}. We give a short argument below for the sake of completeness.
Denote by $f(\mbf{x})$ the largest
$\varepsilon$ for which the open $\varepsilon$-neighborhood of $\mbf x$ is
contained in some $C_i$. Since $f$ is continuous, by compactness, it has a
positive minimum $\delta$. Now $C_i'$ can be the set of all those points whose
open $\delta$-neighborhood is in $C_i$. The condition about antipodal points is
automatically satisfied by $C_i'\subseteq C_i$. Thus the closed set version of
the theorem can be applied to the sets $C_i'$ and the conclusion implies
the same conclusion for the sets $C_i$, again by $C_i'\subseteq C_i$.)
\medskip

The modified version we need is the following.

\medskip
\par
\noindent
{\bf Bacon-Tucker theorem, second form.}
{\it If $C_1,\dots,C_{h+1}$ are open subsets of $S^h$, $$\cup_{i=1}^{h+1}
(C_i\cup (-C_i))=S^h, \ \ \forall i:\ C_i\cap(-C_i)=\emptyset,$$
and $j\in\{0,\dots,h+1\}$, then there is an $\mbf{x}\in S^h$ such that
$$\mbf{x}\in C_i\ {\rm for}\ i\le j\ {\rm and}
\ -\mbf{x}\in C_i\ {\rm for}\ i>j.$$}
\medskip
\par

\proof
Let $D_{h+2}=S^h\setminus\left(\cup_{i=1}^{h+1} C_i\right)$. Then $D_{h+2}\cap
(-D_{h+2})=\emptyset$ by the first condition on the sets $C_i$. Since
$D_{h+2}$ is closed there is some $\varepsilon>0$ bounding the distance
of any pair of points $\mbf x\in D_{h+2}$ and $\mbf y\in -D_{h+2}$ from
below. Let $C_{h+2}$ be the open
${\varepsilon\over 2}$-neighborhood of $D_{h+2}$. Then
the open sets $C_1,\dots,C_{h+2}$ satisfy the conditions (of the open set
version) of the Bacon-Tucker theorem. Therefore its conclusion
holds. Neglecting the set $C_{h+2}$ in this conclusion the proof is
completed for $j>0$. To see the statement for $j=0$ one can take the
negative of the value $\mbf x$ guaranteed for $j=h+1$. \hfill\qed

\medskip
\par
\noindent
{\bf Proof of Theorem~\ref{cklm}.}
Let $G$ be a graph with $\chi(G)=\coind(B_0(G))+1=t$ and fix an arbitrary
proper $t$-coloring $c: V(G)\to [t]$.
Let $g:S^{t-1}\to B_0(G)$ be a $\2$-map that exists by $\coind(B_0(G))=t-1$.

We define
for each color $i\in [t]$ an open set $C_i$ on
$S^{t-1}$. For $\mbf{x}\in S^{t-1}$ we let $\mbf{x}$ be an element of $C_i$
iff the minimal simplex $S_{\msbf{x}}\uplus T_{\msbf{x}}\in B_0(G)$ whose body
contains $g(\mbf{x})$ has a vertex $v\in S_{\msbf{x}}$ for which $c(v)=i$.
These $C_i$'s are open. If an $\mbf{x}\in S^{t-1}$ is not covered by any
$C_i$ then $S_{\msbf{x}}$ must be empty, which also implies
$T_{\msbf{x}}\neq\emptyset$. Since $S_{-\msbf{x}}=T_{\msbf{x}}$ this further
implies $-\mbf x\in\cup_{i=1}^t C_i$, thus $\cup_{i=1}^t (C_i\cup
(-C_i))=S^{t-1}$ follows.

If a set $C_i$ contained an antipodal pair $\mbf{x}$ and $-\mbf{x}$ then
$S_{-\msbf{x}}=T_{\msbf{x}}$ would contain an $i$-colored vertex $u$, while
$S_{\msbf{x}}$ would contain an $i$-colored vertex $v$. Since
$S_{\msbf{x}}\uplus T_{\msbf{x}}$ is a simplex of $B_0(G)$
$u$ and $v$
must be adjacent, contradicting that $c$ is a proper coloring.
Thus the $C_i$'s satisfy the conditions of (the second form of) the
Bacon-Tucker theorem.

Let $j=|A|=l$ and relabel the colors so that colors $1,\dots,j$ be in $A$, and
the others be in $B$. The indices of the $C_i$'s are relabeled accordingly.
We apply
the second form of the Bacon-Tucker theorem with $h=t-1$. It
guarantees the existence of an $\mbf{x}\in S^{t-1}$ with the property
$\mbf{x}\in C_i$ for $i\in A$ and
$-\mbf{x}\in C_i$ for $i\in B$. Then $S_{\msbf{x}}$ contains vertices
$u_1,\dots,u_l$ with $c(u_i)=i$ for all $i\in A$ and
$S_{-\msbf{x}}=T_{\msbf{x}}$ contains vertices $v_1,\dots,v_m$ with
$c(v_i)=l+i$ for all $(l+i)\in B$. Since all vertices of $S_{\msbf{x}}$ are
connected to all vertices in $T_{\msbf{x}}$ by the definition of $B_0(G)$, they
give the required completely multicolored $K_{l,m}$ subgraph. \hfill\qed

\bigskip

\subsection{Graphs that are subject of the colorful $K_{l,m}$ theorem}
\label{gr}

Let us put the statement of Theorem~\ref{cklm} into the perspective of our
earlier work in \cite{ST}. There we investigated the local chromatic number of
graphs that is defined in \cite{EFHKRS} as
$$\psi(G):=\min_c \max_{v\in V(G)} |\{c(u): uv \in E(G)\}|+1,$$
where the minimum is taken over all proper colorings $c$ of $G$. That is
$\psi(G)$ is the minimum number of colors that must appear in the closed
neighborhood of some vertex in any proper coloring.

With similar techniques to those applied in this paper we have shown
using Ky Fan's theorem that $\coind(B_0(G))\ge t-1$ implies that $G$ must
contain a completely multicolored
$K_{\lceil t/2\rceil,\lfloor t/2\rfloor}$ subgraph in {\em any} proper coloring
with the colors alternating with respect to their natural order on the two
sides of this complete bipartite graph. This is the Zig-zag
theorem in \cite{ST} we already referred to in Remark 1 in
Section~\ref{gdol}.
The Zig-zag theorem implies that any graph satisfying its condition must have
$\psi(G)\ge \lceil t/2\rceil +1$. In \cite{ST} we have shown for several
graphs $G$ for which $\chi(G)=\coind(B_0(G))+1=t$ that it can
be colored with $t+1$ colors so that no vertex has more than $\lfloor
t/2\rfloor +2$ colors in
its closed neighborhood. When $t$ is odd, this established the exact value
$\psi(G)=\lceil t/2\rceil+1$ for these graphs.
For odd $t$ this also means that the only
type of $K_{l,m}$ subgraph with $l+m=t$ that must appear completely
multicolored when
using $t+1$ colors is the $K_{\lceil t/2\rceil,\lfloor t/2\rfloor}$ subgraph
guaranteed by the Zig-zag theorem (apart from the empty graph $K_{t,0}$).
If we use only $t$ colors, however, then the
situation is quite different. It is true for any graph $F$ that if it is
properly colored with
$\chi(F)$ colors then each color class must contain a vertex that sees all
other colors in its (open) neighborhood. If it were not so, we could completely
eliminate a color class by recoloring each of its vertices to a color which is
not present in its neighborhood. In the context of local chromaticity
this means that if $\psi(F)<\chi(F)$ then it can only be
attained by a coloring that uses strictly more than $\chi(F)$ colors.
Now Theorem~\ref{cklm} says that if $G$ satisfies $\chi(G)=\coind(B_0(G))+1=t$
then all $t$-colorings give rise not only to completely multicolored $K_{\lceil
t/2\rceil,\lfloor t/2\rfloor}$'s that are guaranteed by the Zig-zag theorem,
and $K_{1,t-1}$'s that must appear in any optimal coloring (with all possible
choices of the color being on the single vertex side), but to all possible
completely multicolored complete bipartite graphs on $t$ vertices.

\bigskip

To conclude this subsection we list explicitly some classes of
graphs $G$ that satisfy the $\chi(G)=\coind(B_0(G))+1$ condition
of Theorem~\ref{cklm}. We recall that in \cite{ST} we used the term
{\em topologically $t$-chromatic} for graphs $G$ with $\coind(B_0(G))\ge t-1$.
With this notation we are listing graphs $G$ that are topologically
$\chi(G)$-chromatic, i.e., for which this specific lower bound on their
chromatic number is tight.

The
subtitles below refer to graph classes that have at least {\it some} members
with the required property. For the precise formulation see
Corollary~\ref{list} below.

\subsubsection{Kneser graphs and Schrijver graphs}

Kneser graphs $\KG(n,k)$ with $t=n-2k+2$ form the probably best known
class of graphs for which the conditions in Theorem~\ref{cklm} hold, cf.,
e.g., \cite{Mat}. Another such class is formed by their vertex color-critical
induced subgraphs called Schrijver graphs $\SG(n,k)$,
also with $t=n-2k+2$. The graph $\SG(n,k)$ is just the general Kneser graph
$\KG(\F)$ for the set system $\F$
consisting of exactly those $k$-subsets of $[n]$ that contain no cyclically
consecutive elements, i.e., neither a pair $\{i,i+1\}$, nor $\{1,n\}$.
Schrijver graphs were introduced by Schrijver in \cite{Schr} and
are discussed in detail in several subsequent papers, here we refer to
\cite{Mat} and \cite{ST} for further information.

The argument presented in Remark~2 proves that the general Kneser graph
$G=\KG(\F)$ also satisfies $\chi(G)=\coind(B_0(G))+1$ as long as
$\chi(\KG(\F))=\cd_2(\F)$.
Note that while the graphs $\KG(n,k)$ are included
in the latter family, the graphs $\SG(n,k)$ are not if $k>1$.

\subsubsection{Mycielski and generalized Mycielski graphs}

Mycielski graphs were introduced by Mycielski \cite{Myc}. These triangle
free graphs are recursively defined starting from $K_2$ and their
chromatic number increases by $1$ at every iteration. It is also true that
applying the Mycielski construction to any graph the clique number will not
change while the chromatic number increases by $1$.

The generalized Mycielski construction first appeared probably in
\cite{Stieb}, cf.\ also \cite{GyJS}. It also reappears in \cite{LLGS, Tar}.
This construction creates from a graph $G$ its
generalized Mycielskian $M_r(G)$, where $r$ denotes the number of ``levels''
in the construction. (The ordinary Mycielski construction is the $r=2$ special
case.) For
a graph $G$ with vertices $v_1,\dots,v_n$ the vertex set of $M_r(G)$ is
$\{v_1^{(p)},\dots,v_n^{(p)}: 0\leq p\leq r-1\}\cup\{z\}$ and the pair
$v_i^{(p)}v_j^{(q)}$ forms an edge iff $v_iv_j\in E(G)$ and either
$p=q=0$ or $p=q\pm 1$. The additional vertex $z$ is connected to
$v_1^{(r)},\dots,v_n^{(r)}$.

When applying this construction to an arbitrary graph, the clique number
does not increase (except in the trivial case when $r=1$) while the chromatic
number may or may not increase. If it does it increases by $1$.
Generalizing Stiebitz's result \cite{Stieb} (see also in \cite{GyJS,Mat})
Csorba \cite{CsMyc} proved that $B(M_r(G))$ is $\2$-homotopy equivalent to
the suspension of $B(G)$ for every graph $G$.
(Csorba's result is in terms of the so-called
homomorphism complex ${\rm Hom}(K_2,G)$ but this is known to be $\2$-homotopy
equivalent to $B(G)$ by results in \cite{Cs, MZ, Ziv}.) Together with
Csorba's already mentioned other result in \cite{Cs} stating the
$\2$-homotopy equivalence of $B_0(G)$ and the suspension of $B(G)$, the
foregoing implies that
if a graph $G$ satisfies $\chi(G)=\coind(B_0(G))+1$, then
the analogous equality will also hold for the graph $M_r(G)$.
In this case the chromatic number does increase by
$1$. Thus iterating the construction $d$ times (perhaps with varying
parameters $r$) we arrive to a graph which has chromatic number $\chi(G)+d$
and still satisfies that its chromatic number is equal to its third (in fact,
for $d>0$ also the fourth) lower bound given in (\ref{eq:chibb0}).
For further explanation of these relations we refer to \cite{ST}.

\subsubsection{Borsuk graphs}

The lower bounds of (\ref{eq:chibb0}) are also tight
for Borsuk graphs $B(n,\alpha)$ with large enough $\alpha<2$
that are defined on $S^{n-1}$ as
vertices and edges are formed by those pairs that are at least distance
$\alpha$ apart, cf.\ \cite{EH,LLgomb}. The chromatic number of these graphs
is $t=n+1$. The paper $\cite{LLgomb}$ shows that
some finite subgraphs of $B(t-1,\alpha)$ also have the required properties.
\bigskip
\bigskip

All the above examples are more or less standard. There are, however, two more,
less known examples the corresponding properties of which are implicit in
\cite{ST}. Below we define these two types of graphs.
Finally, one more family will be given the members of which are
well-known graphs, but their relevance in the present context is probably not
widely known.

\subsubsection{Homomorphism universal graphs for local colorings} \label{Umr}

\begin{defi}\label{defi:Umr} {\rm (\cite{EFHKRS})}
For positive integers $r\leq m$ the graph $U(m,r)$ is defined as follows.
\begin{eqnarray*}
V(U(m,r))&=&\{(i,A): i\in [m], A\subseteq [m], |A|=r-1, i\notin A\}\\
E(U(m,r))&=&\{\{(i,A),(j,B)\}: i\in B, j\in A\}
\end{eqnarray*}
\end{defi}

It is shown in \cite{EFHKRS} that these graphs characterize local colorability
in the following sense: a graph $G$ has an $m$-coloring where no closed
neighborhood of any vertex contains more than $r$ colors if and only if $G$
admits a homomorphism to $U(m,r)$.
As mentioned above, in \cite{ST} we showed for several odd-chromatic graphs
satisfying the
conditions of Theorem~\ref{cklm} that their local chromatic number is $\lceil
t/2\rceil +1$ and it is attained with a coloring using $t+1$ colors. It
follows that for odd $t$ the $t$-chromatic graph $U(t+1, {t+3\over 2})$ also
satisfies the conditions of Theorem~\ref{cklm}. (Indeed, by the functoriality
of $B_0(G)$, $\coind(B_0(F))\ge t-1$ and the existence of a homomorphism $F\to
G$ implies $S^{t-1}\to ||B_0(F)||\to ||B_0(G)||$ and thus $\coind(B_0(G))\ge
t-1$. For further details, cf.\ Remark 3 in \cite{ST}.)

The graphs $U(t+1,{t+2\over 2})$ with $t$ even also belong here.
It is proven in \cite{ST2} that $\coind(B_0(U(t+1,{t+2\over 2}))=t-1$, but the
way of proof is rather different than the previous argument above.
The $t$-colorability of these graphs is also easy to check. We also mention the
result from \cite{ST2} according to which the fourth lower bound on the
chromatic number in (\ref{eq:chibb0}) is not tight for these graphs (while it
is for the graphs of the previous paragraphs). This shows that
Theorem~\ref{cklm} in its present form is somewhat stronger than it would be
with the stronger requirement $\chi(G)=\coind(B(G))+2$ in place of
$\chi(G)=\coind(B_0(G))+1$. We needed the second form of the Bacon-Tucker
theorem for obtaining this stronger form.

We mention that $\chi(U(t+1,\lfloor{{t+3}\over 2}\rfloor))=t$ is a special
case of Theorem 2.6 in \cite{EFHKRS}.

\subsubsection{Homomorphism universal graphs for wide colorings}

\begin{defi} \label{Wst}
Let $H_s$ be the path on the vertices $0,1,2,\dots,s$ ($i$ and $i-1$ connected
for $1\le i\le s$) with a loop at $s$.
We define $W(s,t)$ to be the graph with
\begin{eqnarray*}
V(W(s,t))&=&\{(x_1\dots x_t): \forall i\ x_i\in\{0,1,\dots,s\},\exists!i\
x_i=0,\ \exists j\ x_j=1\}\\
E(W(s,t))&=&\{\{x_1\dots x_t,y_1\dots y_t\}: \forall i\ \{x_i,y_i\}\in
E(H_s)\}.
\end{eqnarray*}
\end{defi}

The graphs $W(2,t)$ are defined in \cite{GyJS} in somewhat different
terms. It is shown there that a graph can be colored properly with $t$ colors
so that the neighborhood of each color class is an independent set if and only
if it admits a homomorphism into $W(2,t)$. The described property is
equivalent to having a $t$ coloring where no walk of length $3$ can connect
vertices of the same color.
Similarly, a graph $F$ admits a homomorphism into $W(s,t)$ if and only if it
can be colored with $t$ colors so that no walk of length $2s-1$ can connect
vertices of the same color. Such colorings are called $s$-wide in \cite{ST}.
Graphs having the mentioned property of $W(s,t)$ are also defined in
\cite{GyJS}, though they are not minimal. The graphs $W(s,t)$ are defined and
shown to be minimal with respect to the above property in \cite{ST}.
The $t$-colorability of $W(s,t)$ is obvious: $c(x_1\dots x_t)=i$ if $x_i=0$
gives a proper coloring. It is
also shown in \cite{ST} that several of the above mentioned $t$-chromatic
graphs (e.g., $B(t-1,\alpha)$ for $\alpha$ close enough to $2$ and $\SG(n,k)$
with $n-2k+2=t$ and $n,k$ large enough with respect to $s$ and $t$) admit a
homomorphism to $W(s,t)$. This implies $\coind(B_0(W(s,t)))\ge t-1$ (with
equality, since $\chi(W(s,t))=t$). Thus the graphs $W(s,t)$ form another
family of graphs for which Theorem~\ref{cklm} applies.

\subsubsection{Rational complete graphs} \label{rat}

Our last example of a graph family satisfying the conditions of
Theorem~\ref{cklm} consists of certain rational (or circular) complete graphs
$K_{p/q}$, as they are called, for example, in \cite{HN}.
The graph $K_{p/q}$ is defined for positive integers $p\ge2q$ on the vertex
set $\{0,\dots,p-1\}$ and $\{i,j\}$
is an edge if and only if $q\leq |i-j|\leq p-q$. The widely investigated
chromatic parameter $\chi_c(G)$, called the circular chromatic number of graph
$G$ (cf.\ \cite{Zhu}, or Section 6.1 in \cite{HN}) can
be defined as the infimum of those values $p/q$ for which $G$ admits a
homomorphism to $K_{p/q}$. It is well known that
$\chi(G)-1<\chi_c(G)\leq\chi(G)$ for every graph $G$. In \cite{LLGS} it is
shown that certain odd-chromatic generalized Mycielski graphs can have their
circular chromatic number arbitrarily close to the above lower bound. Building
on this we showed similar results also for odd chromatic Schrijver graphs and
Borsuk graphs in \cite{ST}. As it is also known that $K_{p/q}$ admits a
homomorphism into $K_{r/s}$ whenever $r/s\ge p/q$ (see, e.g., as Theorem 6.3
in \cite{HN}), the above and the functoriality of $B_0(G)$ together imply that
$\coind(B_0(K_{p/q}))+1=\chi(K_{p/q})=\lceil{p/q}\rceil$ whenever
$\lceil{p/q}\rceil$ is odd.

We remark that the oddness condition is crucial here. It also follows from
results in \cite{ST} (cf.\ also \cite{Meunier} for some special cases) that
the graphs $K_{p/q}$ with
$\lceil{p/q}\rceil$ even and $p/q$ not integral do not satisfy the conditions
of Theorem~\ref{cklm}. Here we state more:
the conclusion of Theorem~\ref{cklm} does not hold for these graphs.
Indeed, let us color the vertex $i$ with the color $\lfloor i/q\rfloor+1$.
This is a
proper coloring with the minimal number $\lceil p/q\rceil$ of colors, but it
does not contain a complete bipartite graph with all the even colors on one
side and all the odd colors on the other.

The remaining case of $p/q$ even and integral is not especially interesting
as $K_{p/q}$ with $p/q$ integral is homomorphic equivalent to the complete
graph on $p/q$ vertices and therefore trivially satisfies the
$\chi(G)=\coind(B_0(G))+1$ condition.

\medskip

Taking the contrapositive in the above observation we obtain a new proof of
Theorem 6 in \cite{ST} the special case of which for Kneser graphs and
Schrijver graphs was independently obtained by Meunier \cite{Meunier}. 

\begin{cor} \label{chic} {\rm (\cite{ST}, cf. also \cite{Meunier})}
If $\coind(B_0(G))$ is odd for a graph $G$, then 
$\chi_c(G)\ge\coind(B_0(G))+1$. 
\end{cor}

\proof
If a graph $G$ has
$\chi_c(G)=p/q$, then $G$ admits a homomorphism to $K_{p/q}$,
thus by the functoriality of $B_0(G)$ we have $\coind(B_0(K_{p/q}))\ge
\coind(B_0(G))$. Then $\lceil p/q\rceil=\chi(K_{p/q})\ge
\coind(B_0(K_{p/q}))+1\ge\coind(B_0(G))+1$. If in addition
$\coind(B_0(G))+1>p/q$, 
then all the previous inequalities hold with equality by the integrality of
the coindex. Then $K_{p/q}$
satisfies the conditions of Theorem~\ref{cklm}, thus it must satisfy its
conclusion. But we just have seen that this is not so if $\lceil p/q\rceil$ is
even and p/q is not integral. Thus $\coind(B_0(G))+1$ cannot be even in this
case.  \hfill\qed

We remark that the consequences of the above result include a partial solution
of two 
conjectures about the circular chromatic number that are mentioned in
\cite{Zhu}. For a detailed discussion of implications and references we refer
to \cite{ST}.  

The proof of Corollary~\ref{chic} in \cite{ST} and also the proof in
\cite{Meunier} relies on Ky Fan's theorem. The above argument shows that Ky
Fan's theorem can be substituted by (the second form of) the Bacon-Tucker
theorem in obtaining this result. Nevertheless, it may be worth noting, that
the missing bipartite graph in the above optimal coloring of an even-chromatic
$K_{p/q}$ is one the presence of which would also be required by the Zig-zag
theorem.

\bigskip
\subsubsection{The entire collection}

Our examples are collected in the following corollary.

\begin{cor} \label{list}
For any proper $t$-coloring of any member of the folowing $t$-chromatic
families of graphs the property described as the conclusion of
Theorem~\ref{cklm} holds.
\begin{description}
\item[(i)] Kneser graphs $\KG(n,k)$ with $t=n-2k+2$,
\item[(ii)] Schrijver graphs $\SG(n,k)$ with $t=n-2k+2$,
\item[(iii)] Borsuk graphs $B(t-1,\alpha)$ with large enough $\alpha<2$ and
some of their finite subgraphs,
\item[(iv)] $U(t+1, \lfloor{t+3\over 2}\rfloor)$, for any $t\ge 2$,
\item[(v)] $W(s,t)$ for every $s\ge 1, t\ge 2$,
\item[(vi)] Rational complete graphs $K_{p/q}$ for $t=\lceil{p/q}\rceil$ odd,
\item[(vii)] The $t$-chromatic graphs obtained by $1\leq d\leq t-2$ iterations
  of the generalized Mycielski construction starting with a $(t-d)$-chromatic
  version of any graph appearing on the list above.
\end{description}
\end{cor}

\proof
The above explanation and references contain the argument implying that
all the families of graphs $G$ above satisfy $t=\chi(G)=\coind(B_0(G))+1$,
thus Theorem~\ref{cklm} is applicable.
\hfill\qed

\subsection{Generalization of G. Spencer and F. Su's result} \label{sect:SpSu}

Recently Gwen Spencer and Francis Edward Su \cite{Spthes,SpSu} found an
interesting consequence of Ky Fan's theorem. They prove that if the Kneser
graph $\KG(n,k)$ is colored optimally, that is, with $t=n-2k+2$ colors, but
otherwise arbitrarily, then the following holds. Given any bipartition of the
color set $[t]$ into partition classes $B_1$ and $B_2$ that are as equal as
possible (i.e., $\Big{|}|B_1|-|B_2|\Big{|}\leq 1$),
there exists a bipartition of the ground set $[n]$ into $E_1$ and $E_2$, such
that, the $k$-subsets of $E_i$ as vertices of $\KG(n,k)$ are all colored with
colors from $B_i$ and every color in $B_i$ does occur (i=1,2).

Theorem~\ref{cklm} implies an analogous statement where no special
requirement is needed about the sizes of $B_1$ and $B_2$. It can also be
obtained by
simply replacing Ky Fan's theorem by the Bacon-Tucker theorem in Spencer and
Su's argument.

\begin{cor} \label{gSpSu}
Let $t=n-2k+2$ and fix an arbitrary proper $t$-coloring $c$ of the Kneser graph
$\KG(n,k)$ with colors from the color set $[t]$. Let $B_1$ and $B_2$ form a
bipartition of $[t]$, i.e., $B_1\cup B_2=[t]$ and $B_1\cap B_2=\emptyset$.
Then there exists a bipartition $(E_1,E_2)$ of $[n]$ such that for $i=1,2$
we have $\{c(v):v\subseteq E_i\}=B_i$.
\end{cor}

\proof
Set $A=B_1$ and $B=B_2$ and consider the complete
bipartite graph Theorem~\ref{cklm} returns for this bipartition
of the color set. Let the vertices on the two sides of this bipartite
graph be $u_1,\dots,u_{|A|}$ and $v_1,\dots,v_{|B|}$. All vertices
$u_i$ and $v_j$ are subsets of $[n]$. Since $u_i$ is adjacent to
$v_j$ for every $i,j$ we have that $E_1':=\cup_{i=1}^{|A|} u_i$
and $E_2':=\cup_{j=1}^{|B|} v_j$ are disjoint. If there are
elements of $[n]$ that are neither in $E_1'$ nor in $E_2'$ then
put each such element into either one of the sets $E_i'$
thus forming the sets $E_1$ and $E_2$. We show that these sets $E_i$
satisfy our requirements. It follows from the construction that
$E_1\cap E_2=\emptyset$ and $E_1\cup E_2=[n]$. It is also clear
that all colors from $B_1$ appear as the color of some $k$-subset
of $E_1$, namely, the $k$-subsets $u_1,\dots,u_{|A|}$ take care
of this condition. Since $E_2$ is disjoint from $E_1$ no
$k$-subset of $E_2$ can be colored by any of the colors from
$B_1$. Thus each $k$-subset of $E_2$ is colored by a color from
$B_2$, and all these colors appear on some $k$-subset of $E_2$ by
the presence of $v_1,\dots,v_{|B|}$. Exchanging the role of $E_1$
and $E_2$ we get that all $k$-subsets of $E_1$ are colored with
some color of $B_1$ and the proof is complete. \hfill\qed

\medskip
\par
\noindent {\it Remark 3.} The same argument proves a similar statement for the
general Kneser graph $\KG(\F)$ in place of $\KG(n,k)$ as long as we have
$t=\chi(\KG(\F))=\coind(B_0(\KG(\F)))+1$. Such graphs include the Schrijver
graphs $SG(n,k)$ with $t=n-2k+2$ and (by the argument presented in Remark 2)
the graphs $\KG(\F)$ with a family $\emptyset\notin\F$ satisfying
$t=\chi(\KG(\F))=\cd_2(\F)$.
\hfill$\Diamond$

\end{document}